\documentclass[12pt,a4paper]{article}
\usepackage[utf8]{inputenc}
\usepackage[T1]{fontenc}
\usepackage[english]{babel}
\usepackage{graphicx}%
\usepackage{multirow}%
\usepackage{amsmath,amssymb,amsfonts}%
\usepackage{amsthm}%
\usepackage{mathrsfs}%
\usepackage[title]{appendix}%
\usepackage{xcolor}%
\usepackage{textcomp}%
\usepackage{manyfoot}%
\usepackage{booktabs}%
\usepackage{algorithm}%
\usepackage{algorithmicx}%
\usepackage{algpseudocode}%
\usepackage{listings}%
\usepackage{mathtools}
\usepackage{amsmath}
\usepackage{amssymb}
\usepackage{hyperref}
\hypersetup{
	colorlinks=true,
	linkcolor=blue,
	filecolor=magenta,      
	urlcolor=blue,citecolor=blue,
	pdftitle={Overleaf Example},
	pdfpagemode=FullScreen,
}
\usepackage[margin=2 cm]{geometry}

\theoremstyle{plain}
\newtheorem{theorem}{Theorem}
\newtheorem{lemma}{Lemma}

\theoremstyle{definition}

\newtheorem{demo}{Proof}

\usepackage{float}

\author{P.~Tiebekabe $^{1,2}$, K. R. Kakanou $^{2}$, H. Ben Yakkou $^{3}$}

\title{Narayana numbers as product of three repdigits in base $g$}

\date{ }

\begin{document}
	
	\maketitle
	\begin{abstract}
		In this paper, we show that there are only finitely many Narayana's 
		numbers which can be written as product of three repdigits in base 
		$g$ with $g \geq  2$. Moreover, for $2 \leq  g \leq  10$, we 
		determine all these numbers. 
	\end{abstract}
	

	\section{Introduction}
	\label{intro}
	
	The problems of the terms of linear recurrence sequences written as a product of 
	repdigits in any base have been intensely studied by several researchers 
	specialized in Number Theory. In this article, we consider the linear recurrent 
	sequence of third order, The Narayana's Cows numbers defined as follows: 
	\[\mathcal{N}_n=\mathcal{N}_{n-1}+\mathcal{N}_{n-3}\quad \text{for}\quad n 
	\geq   3\quad \text{with}\quad \mathcal{N}_0=0, 
	\mathcal{N}_1=\mathcal{N}_2=1.\]	
	
	For more details on the work related to the determination of the terms of linear 
	recurrent sequences which are repdigits in any base, we refer the reader to the 
	following recent results \cite{1}--\cite{4}.   
	
	The concept of Narayana's cows numbers, derived from Indian mythology and 
	Hinduism, holds a significant place in mathematics. These numbers have 
	been extensively studied due to their properties and relationships
	with other mathematical sequences, and their important applications
	in other various fields such as  cryptography, coding theory,
	and graph theory. In this paper, we delve into a fascinating aspect 
	of Narayana numbers by examining their representation as product
	s of three repdigits in base $g$ with $g \geq2$.
	
	Repdigits, which consist of repeated digits, have garnered attention for
	their mathematical properties and patterns. In a fixed base $g \geq 2$,
	a repdigit has the following form, 
	\[ \sum_{i=0}^{n-1} d \times g^{i} = d \times \cfrac{g^n-1}{g-1},\]
	where $1 \leq d \leq g-1$ and $n$ a positive integer.
	
	The proofs of our main results are based on a double application of
	Baker's method and on a reduction algorithm using computations based 
	on continued fractions. The method used to 
	determine Narayana numbers, which are products of three repdigits is
	similar to that used by Ad\'edji \cite{2} and by  Ad\'edji et al.\ \cite{ 3}.
	
	The present paper is organized as follows: in Section \ref{sec2}, we present 
	our main results, Section \ref{sec3} is devoted to reminding necessary 
	results for the proofs of our results, and in Section \ref{sec4},
	we prove our results.

	\section{Statement of main results}\label{sec2}
	In this section, we state all the main results obtained in this paper.
	
	\begin{theorem}\label{MT1}
		Let $g    \geq   2$ be an integer.
		Then the Diophantine equation 
		\begin{equation}\label{Meq1}
			\mathcal{N}_k=d_1\cfrac{g^\ell-1}{g-1}\cdot 
			d_2\cfrac{g^m-1}{g-1}\cdot d_3\cfrac{g^n-1}{g-1}
		\end{equation}
		has only finitely many solutions in integers  $k, d_1, d_2, d_3, \ell
		, m, n$ such that $1 \leq d_i \leq g-1$ for $i=1,2,3$ and 
		$n   \geq m \geq \ell \geq 1$. Further, we have
		\[
		n<5.91 \times 10^{49}\log^9 g\quad\text{and}\quad k<4.73 \times 
		10^{50}\log^{10}g.
		\]
	\end{theorem}
	Under the notation and assumptions of the Theorem \ref{MT1}, if \eqref{Meq1} 
	holds for $\biggl(k,d_1,d_2,d_3,\ell,m,n \biggr) $, then we write 
	\begin{align*}
		\mathcal{N}_{k}=[a,b,c]_g=a\times b \times c,
	\end{align*}
	where
	\[ a=d_1 \times \cfrac{g^{\ell}-1}{g-1}=\overline{\underbrace{ d_1 
			\cdots d_{1_{g}}}_{\ell -1 \ \text{times}}}, \quad b= d_2 \times 
	\cfrac{g^{m}-1}{g-1}=\overline{\underbrace{ d_2 \cdots 
			d_{2_{g}}}_{m -1 \ \text{times}}}, \quad  \text{and}  
	\quad c=d_3 \times \cfrac{g^{n}-1}{g-1}=\overline{\underbrace{ 
			d_3 \cdots d_{3_{g}}}_{n -1 \ \text{times}}}.\]
	In the following theorem, we completely and explicitly give all 
	solutions of the equation \eqref{Meq1} corresponding to $2 \leq g \leq 10.$
	\begin{theorem}\label{MT11}
		The only Narayana numbers which are a product of three repdigits in 
		base $g$ with $2   \leq   g  \leq 10$ are
		\[  \{ 1, 2, 3, 4, 6, 9, 13, 28, 60, 88, 129, 189   \}.\] 
		More precisely, we have 
		\begin{table}[H]
			\caption{Narayana numbers which are a product of three
				repdigits in base $g$, $2 \leq g \leq 10.$}
			\label{Narayana_repdigits}
			\begin{center}
				\begin{tabular}{| c| c | p{12.5cm}| }
					\hline \rule{0pt}{1 \normalbaselineskip}
					$k$ &  $\mathcal{N}_{k}$ & $[a,b,c]_{g}$ 
					\\ 
					\hline \rule{0pt}{1 \normalbaselineskip}
					$1, 2 ,3 $ & $1$ & $ \left[ 1,1,1 \right]_{g} \quad 
					\text{for} \quad g= 2,\ldots, 10.  $  
					\\
					\hline \rule{0pt}{1 \normalbaselineskip}
					$4$ & $2$  & $\left[ 1,1,2 \right]_{g} \quad \text{for} 
					\quad g= 3,\ldots, 10.$
					\\
					\hline \rule{0pt}{1 \normalbaselineskip}
					$5 $ & $ 3 $ & $ \left[ 1,1,11 \right]_{2}$, $\quad  
					\left[1,1,3 \right]_{g} \quad \text{for} \quad g= 4,
					\ldots, 10$. 
					\\
					\hline \rule{0pt}{1 \normalbaselineskip}
					$ 6$ & $4 $ & $  \left[ 1,1,11 \right]_{3}$, 
					$\quad \left[ 1,1,4 \right]_{g} \quad \text{for} \quad
					g= 5,\ldots, 10$, $\quad  \left[ 1,2,2 \right]_{g} 
					\quad \text{for} \quad g= 3,\ldots, 10$. 
					\\
					\hline \rule{0pt}{1 \normalbaselineskip}
					$7 $ & $ 6$ & $ \left[ 1,1,11 \right]_{5} $, $\quad 
					\left[ 1,2,3 \right]_{g}   \quad \text{for} \quad 
					g= 4,\ldots, 10$, $ \quad  \left[ 1,1,6 \right]_{g}   
					\quad \text{for} \quad g= 7,\ldots, 10$.
					\\
					\hline \rule{0pt}{1 \normalbaselineskip}
					$ 8 $ & $ 9$ & $\left[ 1,11,11 \right]_{2}$, $\quad 
					\left[ 1,1,111 \right]_{3}$, $\quad  \left[ 1,1,11 \right]_{8} $,
					$\quad \left[ 1,1,9 \right]_{10}$, $\quad \left[ 1,3,3 \right]_{g}   
					\quad \text{for} \quad g= 4,\ldots, 10 $. 
					\\
					\hline \rule{0pt}{1 \normalbaselineskip}
					$ 9$ & $ 13$ & $\left[ 1,1,111 \right]_{3}	 $ \\
					\hline \rule{0pt}{1 \normalbaselineskip}
					$ 11 $ & $28 $ & $\left[ 1,1,44 \right]_{6}$, $\quad 
					\left[ 1,2,22 \right]_{6}$, $\quad \left[ 1,4,11 \right]_{6}$, 
					$\quad \left[ 2,2,11 \right]_{6}$, $\quad \left[ 1,4,7 \right]_{g} 
					\quad \text{for} \quad g= 8,\ldots, 10$,
					$\left[ 2,2,7 \right]_{g}   \quad \text{for} 
					\quad g= 8,\ldots, 10$.
					\\
					\hline \rule{0pt}{1 \normalbaselineskip}
					$ 13$ & $ 60 $ & $\left[ 2,2,33 \right]_{4}$, 
					$ \quad \left[ 2,3,22 \right]_{4}$, $ \quad 
					\left[ 1,1,66 \right]_{9}$, $ \quad \left[ 1,2,33 \right]_{9}$,
					$\quad\left[ 1,3,22 \right]_{9}$,$ \quad
					\left[ 1,6,11 \right]_{9}$, $ \quad 
					\left[ 2,3,11 \right]_{9}$, $\quad
					\left[ 2,5,6 \right]_{g} \quad \text{for} \quad g= 7,\ldots, 
					10$, $\quad \left[ 3,4,5 \right]_{g} \quad \text{for} 
					\quad g= 6,\ldots, 10$. 
					\\
					\hline \rule{0pt}{1 \normalbaselineskip}
					$ 14$ & $ 88$ & $ \left[ 1,1,88 \right]_{10}$, $\quad
					\left[ 1,2,44 \right]_{10}$, $\quad \left[ 1,4,22 \right]_{10}$,
					$\quad \left[ 1,8,11 \right]_{10}$, $\quad 
					\left[ 2,2,22 \right]_{10}$, $\quad \left[ 2,4,11 \right]_{10}$.
					\\
					\hline \rule{0pt}{1 \normalbaselineskip}
					$ 15 $ & $129 $ & $ \left[ 1,1,333 \right]_{6}$, 
					$\quad \left[ 1,3,111 \right]_{6}$. 
					\\
					\hline \rule{0pt}{1 \normalbaselineskip}
					$16 $ & $189 $ & $ \left[ 1,11,111111 \right]_{2}$, 
					$\quad \left[ 1,3,333 \right]_{4}$, $\quad 
					\left[ 3,3,111 \right]_{4} $, $\quad \left[ 3,3,33 \right]_{6}$,
					$ \quad \left[ 1,3,77 \right]_{8}$, $\quad 
					\left[ 1,7,33 \right]_{8}$, $ \quad \left[ 3,7,11 \right]_{8}$,
					$\quad \left[3,7,9 \right]_{10}$. 
					\\
					\hline 
					
				\end{tabular}
				
			\end{center}
		\end{table}
		
	\end{theorem}
	
	\section{Preliminary Results}\label{sec3}
	In this section, we give some notation and recall certain definitions and 
	results required for the proof of our main results.
	
	\subsection{Some properties of Narayana sequence}
	Narayana's cows sequence comes from a problem with cows proposed by Indian 
	mathematician Narayana in the  14th century. In this problem, we assume that
	there is a cow at the beginning and each cow produces a calf every year from 
	the 4th year. Narayana's cow problem counts the number of calves produced each
	year \cite{5}. 
	
	The characteristic polynomial of Narayana's cows sequence $  \{ \mathcal{N}_n
	\}_{ n \geq0}$ is 
	\[
	\varphi(x)=x^3-x^2-1.
	\]
	Furthermore, the zeros of $\varphi(x)$ are 
	\begin{align*}
		\alpha_{\mathcal{N}} &=  \cfrac{1}{3}\biggl( \sqrt[3]{  
			\cfrac{1}{2}(29-3\sqrt{93})} + \sqrt[3]{  
			\cfrac{1}{2}(3\sqrt{93}+29)}+1 \biggr),
		\\
		\beta_{\mathcal{N}} & =    \cfrac{1}{3}-  \cfrac{1}{6}\biggl( 
		1-i\sqrt{3} \biggr)\sqrt[3]{  \cfrac{1}{2}(29-3\sqrt{93})}- 
		\cfrac{1}{6}\biggl( 1+i\sqrt{3} \biggr)\sqrt[3]{ 
			\cfrac{1}{2}(3\sqrt{93}+29)},
		\\
		\gamma_{\mathcal{N}} &=    \cfrac{1}{3}-  
		\cfrac{1}{6}\biggl( 1+i\sqrt{3} \biggr)\sqrt[3]{ 
			\cfrac{1}{2}(29-3\sqrt{93})}-  \cfrac{1}{6}\biggl( 
		1-i\sqrt{3} \biggr)\sqrt[3]{  
			\cfrac{1}{2}(3\sqrt{93}+29)}.
	\end{align*}
	Then, the Narayana sequence can be obtained by Binet formula
	\begin{equation}\label{eq4}
		\mathcal{N}_n=a_{\mathcal{N}} \alpha_{\mathcal{N}}^n+ b_{\mathcal{N}}
		\beta_{\mathcal{N}}^n + c_{\mathcal{N}} \gamma_{\mathcal{N}}^n.
	\end{equation}
	From  the three initial values of Nayarana sequence, and using Vieta's theorem, 
	one has 
	
	\begin{equation}\label{eq6}
		a_{\mathcal{N}}=  \cfrac{ \alpha_{\mathcal{N}}^2}{ 
			\alpha_{\mathcal{N}}^3+2},\quad	b_{\mathcal{N}}=  
		\cfrac{ \beta_{\mathcal{N}}^2}{ \beta_{\mathcal{N}}^3+2}
		,\quad \text{and}\quad c_{\mathcal{N}}=  \cfrac{ 
			\gamma_{\mathcal{N}}^2}{ \gamma_{\mathcal{N}}^3+2}.
	\end{equation}
	The minimal polynomial of $a_{\mathcal{N}}$ over $\mathbb{Z}$ 
	is $31x^3-3x-1$.
	
	Setting $\Pi(n)=\mathcal{N}_n-a_{\mathcal{N}}
	\alpha_{\mathcal{N}}^n=b_{\mathcal{N}} 
	\beta_{\mathcal{N}}^n+c_{\mathcal{N}} 
	\gamma_{\mathcal{N}}^n$, we notice that
	\begin{equation}\label{Pi_n}
		\biggl|\Pi(n)\biggr| <  
		\cfrac{1}{ \alpha_{\mathcal{N}}^{n/2}} 
		\quad \text{for all}  \ n \geq 1. 
	\end{equation}
	We note that the characteristic polynomial has a 
	real zero $\alpha_{\mathcal{N}}(> 1)$ and two  complex
	conjugate zeros $ \beta_{\mathcal{N}}$ and
	$ \gamma_{\mathcal{N}}$ with $| \beta_{\mathcal{N}}|= 
	| \gamma_{\mathcal{N}}|<1$. In fact, $\alpha_{\mathcal{N}}\approx 1.46557$
	. We also have the following property of $(\mathcal{N}_n)_{n \geq   0}$.
	
	\begin{lemma}\label{lem1}
		For the sequence $(\mathcal{N}_n)_{n \geq   0}$, we have, 
		\[ 
		\alpha_{\mathcal{N}}^{n-2}  \leq    \mathcal{N}_n  \leq   
		\alpha_{\mathcal{N}}^{n-1},\quad\text{for}\quad n  \geq    1.
		\] 
	\end{lemma}
	\begin{demo}
		One can easily prove Lemma ~\ref{lem1} using induction on $n$. 
	\end{demo}  
	
	Let $ \mathbb{K}_{   \varphi}:=\mathbb{Q}(\alpha_{\mathcal{N}}, 
	\beta_{\mathcal{N}}) $ be the splitting field of the polynomial $
	\varphi $ over $ \mathbb{Q} $. Then, $ [\mathbb{K}_{   \varphi}, 
	\mathbb{Q}]=6 $. Furthermore, $ [\mathbb{Q}(\alpha_{\mathcal{N}})
	:\mathbb{Q}]=3 $. The Galois group of $ \mathbb{K} $ over $ \mathbb{Q} $ 
	is given by
	\begin{equation*}
		\mathcal{G}_{   \varphi}:=\text{Gal}(\mathbb{K/Q})\cong \{(1), 
		(\alpha_{\mathcal{N}}\beta_{\mathcal{N}}), 
		(\alpha_{\mathcal{N}}\gamma_{\mathcal{N}}), 
		(\beta_{\mathcal{N}}\gamma_{\mathcal{N}}), 
		(\alpha_{\mathcal{N}}\beta_{\mathcal{N}}\gamma_{\mathcal{N}})
		, (\alpha_{\mathcal{N}}\gamma_{\mathcal{N}}\beta_{\mathcal{N}})\}
		\cong S_3.
	\end{equation*}
	Thus, we identify the automorphisms of $ \mathcal{G}_{   \varphi} $ with 
	the permutations of the zeros of the polynomial $    \varphi $. For example, 
	the permutation $ (\alpha_{\mathcal{N}}\beta_{\mathcal{N}}) $ corresponds to 
	the automorphisms $ \sigma_{   \varphi}: \alpha_{\mathcal{N}} \to 
	\beta_{\mathcal{N}}, ~\beta_{\mathcal{N}} \to \alpha_{\mathcal{N}}, 
	~\gamma_{\mathcal{N}} \to \gamma_{\mathcal{N}} $.

	\subsection{Linear forms in logarithms}
	
	We begin this subsection with a few reminders about the logarithmic height of
	an algebraic number. Let $\eta$ be an algebraic number of degree $d$, $a_0 >0$
	be the leading coefficient of its minimal polynomial over $   \mathbb{Z} $ 
	and let $\eta=\eta^{(1)},\ldots,\eta^{(d)}$ denote its conjugates.
	The quantity defined by  
	\[
	h(\eta)= \cfrac{1}{d}\biggl(\log |a_0|+\sum_{j=1}^{d}\log\max 
	(1,  |\eta^{(j)}   |   ) \biggr)
	\]
	is called the logarithmic height of $\eta$. Some properties of height are as
	follows. For $\eta_1, \eta_2$ algebraic numbers and $m\in    \mathbb{Z} $, 
	we have
	\begin{align*}
		h(\eta_1 \pm \eta_2) &   \leq     h(\eta_1)+ h(\eta_2) +\log2,
		\\
		h(\eta_1\eta_2^{\pm 1}) &   \leq     h(\eta_1) + h(\eta_2),
		\\
		h(\eta_1^m)&=|m|h(\eta_1).
	\end{align*}
	In particular, if $\eta=  p/q \in    \mathbb{Q} $ is a rational number in 
	its reduced form with $q>0$, then $h(\eta)=\log(\max\{|p|,q\})$.
	
	We can now present the famous Matveev's result used in this study.
	Let $\mathbb{L}$ be a real number field of degree $d_{  \mathbb{L} }$, 
	$\eta_1,\ldots,\eta_s \in \mathbb{L}$ and $b_1,\ldots,b_s \in    \mathbb{Z}  
	\setminus\{0\}$. Let $B    \geq     \max\{|b_1|,\ldots,|b_s|\}$ and
	\[
	\Lambda=\eta_1^{b_1}\cdots\eta_s^{b_s}-1.
	\]
	Let $A_1,\ldots,A_s$ be real numbers such that 
	\[
	A_i    \geq     \max\{d_{  \mathbb{L} } h(\eta_i), |\log\eta_i|, 0.16\},
	\quad i=1,\ldots,s.
	\]
	With the above notation, Matveev \cite{6} proved the following result.
	\begin{theorem}\label{Matveev}
		Assume that $\Lambda\neq 0$. Then
		\[\log|\Lambda|>-1.4\cdot30^{s+3}\cdot s^{4.5}\cdot d_{\mathbb{L}}^2 
		\cdot(1+\log d_{  \mathbb{L} })\cdot(1 +\log B)\cdot A_1\cdots A_s.\]
	\end{theorem}
	We also need the following result from Sanchez and Luca \cite{7}.
	\begin{lemma}\label{lem3}
		Let $r  \geq    1$ and $H>0$ be such that $H>(4r^2)^r$ and 
		$H>L/(\log L)^r$. Then 
		\[ 
		L<2^rH(\log H)^r.
		\]
	\end{lemma}
	\subsection{Reduction method}
	The bounds on the variables obtained via Baker’s theory \cite{8} are too large
	for any computational purposes. To reduce the bounds, we use the reduction 
	method due to Dujella and Peth\H{o} \cite[Lemma 5a]{9}.
	For a real number $X$, $  |X  | :=\min\{  | X-n   | :n\in 
	\mathbb{Z}\}$ stands for the distance of $X$ to the nearest integer.
	
	\begin{lemma}\label{Dujella}
		Let $M$ be a positive integer, $p/q$ be a convergent of the continued 
		fraction expansion of an irrational number $\tau$ such that $q>6M$, 
		and $A, B, \mu$ be some real numbers with $A>0$ and $B>1$. Furthermore,
		let
		\[  
		\varepsilon:=  |\mu q   |-M\cdot   |\tau q  | .
		\]
		If $\varepsilon>0$, then there is no solution to the inequality
		\begin{equation}
			0<|u\tau-v+\mu|<AB^{-w}
		\end{equation}
		in positive integers $u,v$ and $w$ with
		\[ 
		u   \leq    M\; \mbox{and}\; w  \geq   
		\cfrac{\log(Aq/\varepsilon)}{\log B}.
		\]
	\end{lemma}
	\section{Proofs of main results} \label{sec4}
	
	\subsection{Proof of Theorem \ref{MT1}}
	To prove Theorem  \ref{MT1},  we  will use the following lemma  which provides a
	relation  on the size of $k$ versus $n$ and $g$.
	\begin{lemma}\label{ML1}
		
		All solutions of the Diophantine equation \eqref{Meq1} satisfy
		\begin{align*}
			k<8 n\log g.	
		\end{align*}	
	\end{lemma}
	\begin{demo}
		From \eqref{lem1}, we have
		\[   
		\alpha_{\mathcal{N}}  ^{k-2}  \leq  \mathcal{N}_k=d_1\cfrac{g^\ell-1}
		{g-1}\cdot d_2\cfrac{g^m-1}{g-1}\cdot d_3\cfrac{g^n-1}{g-1}   
		\leq (g^n-1)^3<g^{3n}.
		\]
		Taking logarithm on both sides, we get $(k-2)\log  
		\alpha_{\mathcal{N}}  <3n\log g$. Since $n   \geq  2$ and $g   \geq
		2$, we obtain the desired inequality. This ends the proof. 
	\end{demo}
	\begin{demo}[Proof of Theorem \ref{MT1}]
		If $n=1$, then $\ell =m=1$. So, equation \eqref{Meq1} becomes
		\[ \mathcal{N}_k = d_1 d_2 d_3\]
		which implies 
		\[\alpha_{\mathcal{N}}^{k-2}   \leq \biggl( g-1 \biggr)^3 \]
		which leads to 
		\[ k < 2+ 3 \cfrac{\log g}{\log \alpha_{\mathcal{N}}} .\]
		Now, suppose $n   \geq 2$.
		From \eqref{Meq1} and \eqref{eq4}, we have
		\[	\mathcal{N}_k=a_{\mathcal{N}} \alpha_{\mathcal{N}}^k + 
		b_{\mathcal{N}} \beta_{\mathcal{N}}^k+ c_{\mathcal{N}} 
		\gamma_{\mathcal{N}}^k=d_1\cfrac{g^\ell-1}{g-1}\cdot 
		d_2\cfrac{g^m-1}{g-1}\cdot d_3\cfrac{g^n-1}{g-1}.\]
		Which implies
		\begin{equation}\label{1MEq1}
			\begin{aligned}
				a_{\mathcal{N}} \alpha_{\mathcal{N}}^k - 
				\cfrac{d_1d_2d_3 g^{ \ell +m+n}}{(g-1)^3} 
				&= - \cfrac{d_1d_2d_3 \biggl(g^{ \ell +m}  
					+ g^{ \ell +n}+g^{m+n}\biggr)}{(g-1)^3}
				\\
				&  + \cfrac{d_1d_2d_3 \biggl(g^{l}+ g^{m}+ 
					g^{n}\biggr)}{(g-1)^3} - \cfrac{d_1d_2d_3 }
				{(g-1)^3} - \Pi(k).
			\end{aligned}
		\end{equation}	
		Taking the absolute values of both sides of \eqref{1MEq1} and using 
		\eqref{Pi_n}, we get
		
		\begin{equation}\label{1MEq2}
			\begin{aligned}
				\biggl|a_{\mathcal{N}} \alpha_{\mathcal{N}}^k - 
				\cfrac{d_1d_2d_3 g^{\ell+m+n}}{(g-1)^3}\biggr| & <  
				\cfrac{d_1d_2d_3 \biggl(g^{\ell+m}  + g^{\ell+n}+
					g^{m+n}\biggr)}{(g-1)^3} 
				\\
				&  + \cfrac{d_1d_2d_3 \biggl(g^{l}+ g^{m}+ 
					g^{n}\biggr)}{(g-1)^3} + \cfrac{d_1d_2d_3 }
				{(g-1)^3} + \cfrac{1}{\alpha_{\mathcal{N}}^{k/2}} .	 
			\end{aligned}
		\end{equation}	
		Multiplying both sides of \eqref{1MEq2} by $\cfrac{(g-1)^3 }
		{d_1d_2d_3g^{\ell+n+m}}$ and noticing the fact that $\ell  
		\leq m  \leq n$, we get the inequality
		\begin{align*}
			\biggl|\cfrac{(g-1)^3\cdot  a_\mathcal{N} 
				\alpha_{\mathcal{N}}^k\cdot g^{-(\ell+n+m)}}
			{d_1d_2d_3}-1\biggr| & <\cfrac{1}{g^\ell}+
			\cfrac{1}{g^m}+\cfrac{1}{g^n}+\cfrac{1}
			{g^{\ell+m}} +\cfrac{1}{g^{\ell+n}}
			\\
			&+\cfrac{1}{g^{m+n}}+  \cfrac{1}{g^{\ell+m+n}}+ 
			\cfrac{(g-1)^3}{\alpha_{\mathcal{N}}^{k/2} 
				d_1d_2d_3g^{\ell+n+m}}
			\\
			&< 8\cdot g^{-\ell}.
		\end{align*}
		So, we get 
		\begin{align}\label{1MEq3}
			\biggl|\cfrac{a_\mathcal{N} (g-1)^3}{d_1d_2d_3}\cdot
			\alpha_{\mathcal{N}}^k \cdot g^{-(\ell+n+m)}-1\biggr|
			< 8\cdot g^{-\ell}.
		\end{align}
		We put 
		\[ \Gamma_1 := \cfrac{a_\mathcal{N} (g-1)^3}{d_1d_2d_3}\cdot 
		\alpha_{\mathcal{N}}^k \cdot g^{-(\ell+n+m)}-1  .\]
		Let us show $ \Gamma_1  \ne 0$. We proceed by the contrary. 
		Assume that $ \Gamma_1  =0$. Then
		\[ a_\mathcal{N} \alpha_{\mathcal{N}}^k =  \cfrac{d_1d_2d_3}
		{(g-1)^3} \cdot g^{\ell +m+n}\]
		which implies 
		\[ \sigma_{\varphi}\biggl( a_\mathcal{N} \alpha_{\mathcal{N}}^k\
		\biggr)= b_\mathcal{N} \beta_{\mathcal{N}}^k = \cfrac{d_1d_2d_3}
		{(g-1)^3} \cdot g^{\ell +m+n}.\]
		Taking the absolute value, we get 
		\[\biggl|b_\mathcal{N} \beta_{\mathcal{N}}^k\biggr| = \biggl| 
		\cfrac{d_1d_2d_3}{(g-1)^3} \cdot g^{\ell +m+n}\biggr|.\]
		We have $ \biggl|b_\mathcal{N} \beta_{\mathcal{N}}^k  \biggr|
		< 1$ instead $ \biggl|\cfrac{d_1d_2d_3}{(g-1)^3} \cdot 
		g^{\ell +m+n}   \biggr|  > 1$ since $1   \leq \ell  
		\leq m   \leq n$, which leads to a contradiction. 
		Hence $ \Gamma_1  \ne 0$. 
		
		In order to apply Matveev's result to $ \Gamma_1 $, set 
		\begin{align*} 
			&t:=3, \quad \eta _1:=\cfrac{a_\mathcal{N} (g-1)^3}
			{d_1d_2d_3} , \quad \eta _2:=\alpha_{\mathcal{N}} , 
			\quad \eta _3:=g,
			\\
			&\quad b_1:=1, \quad b_2:=k, \quad b_3:=-(\ell +m+n),
		\end{align*}
		and ${\mathbb {K}}:={\mathbb {Q}}(\eta _1,\eta _2, \eta _3)=
		{\mathbb {Q}}(\alpha_{\mathcal{N}} )$ which is a real number 
		field of degree $d_{\mathbb {K}}=3$.
		
		Using properties of the logarithmic height, we get  \[h 
		(\eta_2   )=h  ( \alpha_{\mathcal{N}}  )
		= \cfrac{\log \alpha_{\mathcal{N}}}{3} , \quad h 
		(\eta_3   )=h   (g   )= \log g\] and 
		\begin{align*}
			h(\eta_1) & =h \biggl(\cfrac{a_\mathcal{N} 
				(g-1)^3}{d_1d_2d_3}\biggr)
			\\
			&   \leq h(a_\mathcal{N} ) + h\biggl(\cfrac{ (g-1)^3}{d_1d_2d_3}
			\biggr)\\
			&   \leq \cfrac{1}{3} \log 23 + \log \biggl( \max \biggl\{   (g-1 
			)^3, d_1d_2d_3 \biggr\} \biggr)
			\\
			& < 3 \log (g)+2< 6 \log g \quad \text{since} \ g   \geq 2 . 
		\end{align*}
		Thus, we can take 
		\[ A_1= 18 \log   (g  ) , \quad A_2= \log \alpha_{\mathcal{N}},
		\quad \text{and} \ A_3= 3 \log g.\]
		By Lemma \ref{ML1}, we have $k < 8 n \log  g $, so we put $B= 8 n \log g$.
		\\
		Using Theorem \ref{Matveev}, we see that 
		\begin{align*}
			\log   |\Gamma_1   | > & -  1.4 
			\times 30^6 \times 3^{4.5}\times 3^2   
			(1 + \log 3)(1+\log (8n \log g )) 
			\\
			& \qquad \times (18\log  (g   ))
			(3\log g  \log \alpha_{\mathcal{N}})
			\\
			& > -5.6 \times 10^{13} (1+\log (8n \log g ))(\log^2 g ). 
		\end{align*}
		Comparing the above inequality with \eqref{1MEq3}, we obtain that
		\begin{align*}
			\ell \log g - \log 8 < 5.6 \times 10^{13} (1+\log 
			(8n \log g ))(\log^2 g ).
		\end{align*}
		Since $g   \geq 2$ and $n   \geq 2$, we have
		\[ 1+\log (8n \log g )< 8 \log (n \log g)     \]
		so we get 
		\[\ell < 4.5 \times 10^{14} \log n \log^2 g.\]
		Rewrite \eqref{Meq1}, we get
		\[	\cfrac{ a_\mathcal{N} \alpha_{\mathcal{N}}^k 
			(g-1)}{d_1(g^\ell-1)}+\cfrac{\Pi(k) (g-1)}{d_1(g^\ell-1)}
		=\cfrac{d_2d_3}{(g-1)^2}\biggl(g^{n+m}-g^m-g^n+1 \biggr),\]
		which implies 
		\begin{align}\label{1Meq4}
			\cfrac{ a_\mathcal{N} \alpha_{\mathcal{N}}^k (g-1)}
			{d_1(g^\ell-1)}-  \cfrac{d_2d_3g^{n+m}}{(g-1)^2}=-
			\cfrac{\Pi(k) (g-1)}{d_1(g^\ell-1)}- 
			\cfrac{d_2d_3g^m}{(g-1)^2}-   \cfrac{d_2d_3g^n}
			{(g-1)^2} +  \cfrac{d_2d_3}{(g-1)^2}.
		\end{align}
		Taking the absolute values of both sides of \eqref{1Meq4}, we have
		\[  \biggl|\cfrac{ a_\mathcal{N} \alpha_{\mathcal{N}}^k 
			(g-1)}{d_1(g^\ell-1)}-  \cfrac{d_2d_3g^{n+m}}{(g-1)^2}\biggr|
		< \cfrac{ (g-1)}{d_1(g^\ell-1) \alpha_{\mathcal{N}}^{k/2}} 
		+   \cfrac{d_2d_3g^m}{(g-1)^2}+   \cfrac{d_2d_3g^n}{(g-1)^2}
		+  \cfrac{d_2d_3}{(g-1)^2}. \]
		Dividing both sides of the inequality above by $  
		\cfrac{d_2d_3g^{n+m}}{(g-1)^2}$ and using the fact that 
		$n  \geq 2$, we see that
		\begin{align*}
			\biggl|  \cfrac{(g-1)^3}{d_1d_2d_3(g^\ell-1)}\cdot 
			a_\mathcal{N} \alpha_{\mathcal{N}}^k \cdot 
			g^{-(n+m)}-1\biggr|&  \leq  \cfrac{(g-1)^3}{d_1d_2d_3(g^\ell-1) 
				\alpha_{\mathcal{N}}^{k/2} g^{n+m} }+  \cfrac{1}{g^n}+  
			\cfrac{1}{g^m}+  \cfrac{1}{g^{n+m}}
			\\
			& <4\cdot g^{-m}.
		\end{align*}
		Then, we have
		\begin{align}\label{1MEq5}
			\biggl|\cfrac{(g-1)^3}{d_1d_2d_3(g^\ell-1)}\cdot 
			a_\mathcal{N} \alpha_{\mathcal{N}}^k \cdot g^{-(n+m)}-1
			\biggr|  < \cfrac{4}{g^m}.
		\end{align}
		We put
		\[\Gamma_2=   \cfrac{(g-1)^3}{d_1d_2d_3(g^\ell-1)}\cdot 
		a_\mathcal{N} \alpha_{\mathcal{N}}^k \cdot g^{-(n+m)}-1.\]
		One can check that $\Gamma_2 \ne 0$, proceeding as we do 
		for $\Gamma_1$. Let us apply Matveev's result for $\Gamma_2$.
		Let
		\begin{align*} 
			& t:=3, \quad \eta _1:=  
			\cfrac{(g-1)^3}{d_1d_2d_3(g^\ell-1)}\cdot a_\mathcal{N}
			, \quad \eta _2:=\alpha_{\mathcal{N}} , \quad \eta _3:=g,
			\\
			& \quad b_1:=1, \quad b_2:=k, \quad b_3:=-(m+n),
		\end{align*}
		and ${\mathbb {K}}:={\mathbb {Q}}(\eta _1,\eta _2, \eta _3)=
		{\mathbb {Q}}(\alpha_{\mathcal{N}} )$ of degree $d_{\mathbb {K}}=3$. 
		By Lemma \ref{ML1}, we have $k < 8 n \log  g $, so we put $B= 8 n 
		\log g$. We have \[h   (\eta_2   )=h  ( 
		\alpha_{\mathcal{N}}   )= \cfrac{\log \alpha_{\mathcal{N}}}{3} 
		, \quad h   (\eta_3   )=h   (g   )= \log g,\] 
		and 
		\begin{align*}
			h(\eta_1) & =h \biggl(  \cfrac{(g-1)^3}{d_1d_2d_3(g^\ell-1)}
			\cdot a_\mathcal{N} \biggr)
			\\
			&   \leq h(a_\mathcal{N} ) + h\biggl(  \cfrac{(g-1)^3}
			{d_1d_2d_3(g^\ell-1)} \biggr)
			\\
			&   \leq \cfrac{1}{3} \log 23 + \log \biggl(\max 
			\biggl\{   (g-1   )^3, d_1d_2d_3 \biggr\} \biggr) + 
			h\biggl(\cfrac{1}{g^{\ell}-1} \biggr)
			\\
			&  < 2+ 3 \log   (g-1   ) + \log
			(g^{\ell}-1   )
			\\
			& <   (3+\ell   ) \log (g)+2
			\\
			&  	< (6+\ell) \log g \quad \text{since} \ g   \geq 2 . 
		\end{align*}
		Thus, we can take 
		\[ A_1= (18+3 \ell) \log   (g   ) , \quad A_2= \log
		\alpha_{\mathcal{N}}, \quad \text{and} \ A_3= 3 \log g.\]
		Using Theorem \ref{Matveev}, we see that
		\begin{align*}
			\log    |\Gamma_2    |  > & -  1.4 \times 30^6
			\times 3^{4.5}\times 3^2     (1 + \log 3)(1+\log (8n \log g ))
			\\
			& \qquad \times ((18+3 \ell) \log   (g   ) )
			(3\log g  \log \alpha_{\mathcal{N}})
			\\
			& > -3.1 \times 10^{12} (1+\log (8n \log g ))(\log^2 g )
			(18+3 \ell).
		\end{align*}
		Comparing with \eqref{1MEq5}, we get 
		\[m \log g - \log 4 < 3.1 \times 10^{12} (1+\log (8n \log g ))
		(\log^2 g ) (18+3 \ell). \]
		We have \[1+ \log (8n \log g) <  8\log n \log g \quad \text{and} 
		\quad \ell < 4.5 \times 10^{14} \log n \log^2 g.\]
		So,
		\[m < 3.8 \times 10^{28} \log^2 n \log^4 g .     \]
		Reorganizing \eqref{Meq1}, we get
		\[	  \cfrac{d_3 g^n}{g-1}-  \cfrac{(g-1)^2 \cdot a_\mathcal{N}
			\alpha_\mathcal{N}^k}{d_1d_2(g^\ell-1)(g^m-1)}=  \cfrac{d_3}{g-1}+ 
		\cfrac{\Pi(k)(g-1)^2}{d_1d_2(g^\ell-1)(g^m-1)}.\]
		We have
		\[  \biggl|\cfrac{d_3 g^n}{g-1}-  \cfrac{(g-1)^2 \cdot a_\mathcal{N}
			\alpha_\mathcal{N}^k}{d_1d_2(g^\ell-1)(g^m-1)}  \biggr|  < 
		\cfrac{d_3}{g-1}+  \cfrac{(g-1)^2}{\alpha_{\mathcal{N}}^{k/2}d_1d_2
			(g^\ell-1)(g^m-1)} \]
		by taking the absolute values of both sides of \eqref{1Meq4}.	
		Dividing both sides of the above inequality by $  (d_3g^{n})/(g-1)$
		and using the fact that $n  \geq 2$, we see that
		\begin{align*}
			\biggl|1-   \cfrac{a_\mathcal{N}(g-1)^3}{d_1d_2d_3(g^\ell-1)
				(g^m-1)}\cdot g^{-n}\cdot \alpha_\mathcal{N}^k\biggr| & < 
			\cfrac{1}{g^n}+  \cfrac{1}{g^{n-1}}<  \cfrac{2}{g^{n-1}}.
			\\
			& < 2 \cdot g^{1-n}.
		\end{align*}
		Then, we have
		\begin{align}\label{1MEq6}
			\biggl|\cfrac{a_\mathcal{N}(g-1)^3}{d_1d_2d_3(g^\ell-1)
				(g^m-1)}\cdot g^{-n}\cdot \alpha_\mathcal{N}^k-1\biggr| 
			< 2 \cdot g^{1-n}.
		\end{align}
		We put
		\[\Gamma_3=\cfrac{a_\mathcal{N}(g-1)^3}{d_1d_2d_3(g^\ell-1)(g^m-1)}
		\cdot g^{-n}\cdot \alpha_\mathcal{N}^k-1 .\]
		One can verify that $\Gamma_3 \ne 0$. Let us analyze  Matveev's 
		result for $ \Gamma_3 $. Let  
		\begin{align*} 
			& t:=3, \quad \eta _1:=  \cfrac{a_\mathcal{N}(g-1)^3}
			{d_1d_2d_3(g^\ell-1)(g^m-1)}, \quad \eta _2:=\alpha_{\mathcal{N}} 
			, \quad \eta _3:=g, 
			\\
			&\quad b_1:=1, \quad b_2:=k, \quad b_3:=-n,
		\end{align*}
		and ${\mathbb {K}}:={\mathbb {Q}}(\eta _1,\eta _2, \eta _3)={\mathbb {Q}}
		(\alpha_{\mathcal{N}} )$ of degree $d_{\mathbb {K}}=3$. By Lemma \ref{ML1}, 
		we have $k < 8 n \log  g $, so we put $B= 8 n \log g$. We have 
		\[h   (\eta_2   ) =h  ( \alpha_{\mathcal{N}}   )= 
		\cfrac{\log \alpha_{\mathcal{N}}}{3} , \quad h   (\eta_3   ) 
		=h   (g   )= \log g,\] and 
		\begin{align*}
			h(\eta_1) & =h \biggl(\cfrac{a_\mathcal{N}(g-1)^3}{d_1d_2d_3
				(g^\ell-1)(g^m-1)} \biggr)
			\\
			&   \leq h(a_\mathcal{N} ) + h\biggl(  \cfrac{(g-1)^3}
			{d_1d_2d_3(g^\ell-1)(g^m-1)} \biggr)
			\\
			&   \leq \cfrac{1}{3} \log 23 + \log \biggl(\max 
			\{   (g-1   )^3, d_1d_2d_3 \} \biggr) + 
			h\biggl(\cfrac{1}{g^{\ell}-1} \biggr)+
			h\biggl(\cfrac{1}{g^{m}-1} \biggr)
			\\
			&  < 2+ 3 \log \biggl(g-1 \biggr) + \log \biggl(g^{\ell}-1
			\biggr) +\log \biggl(g^{m}-1 \biggr)
			\\
			& <   (3+\ell+m   ) \log (g)+2
			\\
			&  	< (6+\ell+m) \log g \quad \text{since} \ g   \geq 2 . 
		\end{align*}
		Thus, we can take 
		\[ A_1= 3(6+\ell+m) \log \biggl(g \biggr)  , \quad A_2= \log 
		\alpha_{\mathcal{N}}, \quad \text{and} \ A_3= 3 \log g.\]
		Using Theorem \ref{Matveev}, we see that
		\begin{align*}
			\log  \biggl|\Gamma_3  \biggr|   &> -  1.4 \times 30^6 \times 
			3^{4.5}\times 3^3    (1 + \log 3)(1+\log (8n \log g ))
			\\
			& \qquad \times ((6+ \ell+m) \log \biggl(g \biggr) )(3\log g 
			\log \alpha_{\mathcal{N}})
			\\
			& > -9.31 \times 10^{12} (1+\log (8n \log g ))(\log^2 g ) 
			(6+ \ell+m).
		\end{align*}
		Comparing with \eqref{1MEq6}, we get 
		\[(n-1) \log g - \log 2 < 9.31 \times 10^{12} (1+\log (8n \log g ))
		(\log^2 g ) (6+ \ell+m) \]
		We have \begin{align*}
			& 1+ \log (8n \log g) <  8 \log n \log g, \quad m <  
			3.8 \times 10^{28} \times \log^2 n \log^4 g \\
			& \text{and} \quad \ell < 4.5 \times 10^{14} \log n \log^2 g.
		\end{align*}
		Thus 
		\begin{align*}
			6+\ell+m & < 4.51 \times 10^{14} \log n \log^2 g+  3.8 \times
			10^{28} \times \log^2 n \log^4 g
			\\
			& < 4 \times 10^{28} \times \log^2 n \log^4 g.
		\end{align*}
		So, we have
		\[n < 3 \times 10^{42}  \log^3 n \log^6 g  .\]
		Now we apply Lemma \ref{lem3}, by setting 
		\[ 	r:=3,\quad L:=n\quad \text{and}\quad H:=3\times 10^{42}
		\cdot\log^6 g,\]
		we get
		\begin{align*}
			n&<2^3\cdot 3 \times 10^{42}\cdot\log^6 g\times \log^3 
			\biggl(3\times 10^{42}\cdot\log^6 g \biggr)
			\\
			&< 2.4\times 10^{43}\cdot\log^6 g\cdot (95.6+6\log\log g)^3
			\\
			&<5.91\times 10^{49}\log^9 g. 
		\end{align*}
		Notice that we have used the inequality  $95.6+6\log\log g< 135 
		\log g$ which holds since $g   \geq 2$.
	\end{demo}
	\subsection{Proof of Theorem \ref{MT11}}
	
	Since $2   \leq g   \leq 10$, according to Theorem \ref{MT1}, we have 
	
	\[		\ell   \leq m   \leq n< 1.08 \times 10^{53} \quad \text{and}
	\quad k < 1.99 \times 10^{54} .
	\]
	Consequently, the next step is to reduce the upper bounds above in order to
	identify the set of the interval in which the possible solutions of \eqref{Meq1}
	lie. To do this, we proceed in three steps.
	\subsubsection*{
		Step~1.} 
	
	Using \eqref{1MEq3}, let 
	$$
	\Lambda_1:=-\log(\Gamma_1+1)=(\ell+m+n)\log g -k \log \alpha_{\mathcal{N}}-\log
	\biggl(  \cfrac{( g -1)^3 a_{\mathcal{N}}}{d_1d_2d_3}\biggr).
	$$
	Notice that \eqref{1MEq3} can be rewritten as 
	$$
	\biggl|  e^{-\Lambda_1}-1 \biggr|< \cfrac{8}{ g ^{\ell}}.
	$$
	Observe that $\Lambda_1\neq 0$, since $ e^{-\Lambda_1}-1=\Gamma_1\neq 0$.
	Assume that $\ell  \geq 5$. Then 
	\[
	\biggl|  e^{-\Lambda_1}-1 \biggr|< \cfrac{8}{ g ^{\ell}}<\cfrac{1}{2}.
	\]
	Since $ \biggl|x  \biggr| < 2 \biggl|\mathrm{e}^x -1  \biggr| $, 
	if $ \biggl|x  \biggr| < \cfrac{1}{2}$ holds, then
	$$
	\biggl|  \Lambda_1 \biggr|<\cfrac{16}{ g ^{\ell}}.
	$$
	Substituting $\Lambda_1$ in the above inequality with its value and
	dividing through by $\log \alpha_{\mathcal{N}}$, we get
	$$
	\biggl| (\ell+m+n)\biggl( \cfrac{\log  g }{\log \alpha_{\mathcal{N}}}\biggr)-k+
	\cfrac{\log \biggl(\cfrac{( g -1)^3 a_{\mathcal{N}}}{d_1d_2d_3} \biggr)}
	{\log \alpha_{\mathcal{N}}} \biggr|<\cfrac{16}{ \log \alpha_{\mathcal{N}} 
		g ^{\ell} }.
	$$ 
	Then, we can apply Lemma \ref{Dujella} with the data 
	\begin{align*}
		\tau:=\cfrac{\log  g }{\log \alpha_{\mathcal{N}}},\quad
		\mu:=\cfrac{\log \biggl(\cfrac{( g -1)^3 a_{\mathcal{N}}}{d_1d_2d_3} 
			\biggr) }{\log \alpha_{\mathcal{N}}},\quad A:=\cfrac{16}{ \log 
			\alpha_{\mathcal{N}}},\quad 
		\\
		B:= g ,\quad \text{and}\quad w:=\ell \quad\text{with}\quad 1
		\leq d_1  \leq d_2  \leq d_3 \leq  g -1.
	\end{align*}
	
	We can take $M:=1.99 \times 10^{54}$, since $k < 8 n \log g <1.99\times 
	10^{54}$. So, for the remaining proof, we use \textit{Mathematica} to apply
	Lemma~\ref{Dujella}. For the computations, if the first convergent $q_t$ 
	is such that $q_t > 6M$ does not satisfy the condition $\varepsilon>0$, 
	then we use the next convergent until we find the one that satisfies the 
	conditions. Thus, we have that 
	
	\begin{table}[ht]
		\caption{Upper bound on $\ell$}
		\label{Bound_ell}
		\begin{center}
			\begin{tabular}{|c||c|c|c|c|c|c|c|c|c|}
				
				\hline \rule{0pt}{1.5 \normalbaselineskip} 
				$ g $ &   {2} &   {3} &   {4} &   {5} &   {6} & 
				{7} &   {8} &   {9} &   {10} 
				\\
				\hline 
				\hline \rule{0pt}{1.5 \normalbaselineskip}
				$q_t$  & $q_{118}$ & $q_{100}$ & $q_{110}$ & $q_{115}$ &
				$q_{90}$ & $q_{106}$ & $q_{112}$ & $q_{102}$ 
				& $q_{96}$
				\\
				\hline \rule{0pt}{1.5 \normalbaselineskip}
				$\varepsilon  \geq$ &0.36 & 0.26 & 0.03 & 0.01&
				0.06 & 0.001 & 0.0019 & 0.005 &0.01
				\\ 
				
				\hline \rule{0pt}{1.5 \normalbaselineskip}
				$\ell  \leq$  & 194  & 121   & 99   & 87       
				& 76    & 72   & 67      &   62 & 59   
				\\ 
				\hline
			\end{tabular}
			
		\end{center}
	\end{table}
	
	\vspace{3mm}
	Therefore
	$$
	1 \leq \ell  \leq \cfrac{\log((16/\log \alpha_{\mathcal{N}}) \cdot 
		q_{118}/0.36)}{\log 2} \leq 194.
	$$ 
	\subsubsection*{Step~2.} 
	In this step, we have to reduce the upper bound on $m$. To do this, 
	let us consider
	$$
	\Lambda_2:=-\log (\Gamma_2+1)=(m+n)\log g - k \log \alpha_{\mathcal{N}} 
	+\log \biggl(\cfrac{ ( g -1)^3 a_{\mathcal{N}}}{d_1d_2d_3 \cdot 
		( g^\ell-1   )}  \biggr).
	$$
	Thus inequality \eqref{1MEq5} becomes 
	$$
	\biggl|  e^{-\Lambda_2}-1 \biggr|< \cfrac{4}{ g ^{m}}<\cfrac{1}{2},
	$$
	which holds for $m \geq 4$. It follows that 
	\begin{align}\label{Duj}
		\biggl|(m+n)\cfrac{\log g}{\log \alpha_{\mathcal{N}}} - k  
		+\cfrac{\log\biggl(\cfrac{ ( g -1)^3 
				a_{\mathcal{N}}}{d_1d_2d_3 \cdot   ( 
				g^\ell-1   ) }\biggr)}{\log 
			\alpha_{\mathcal{N}}} \biggr|< \cfrac{8}{ g ^{m} \log
			\alpha_{\mathcal{N}} }.
	\end{align}
	So, the conditions of Lemma \ref{Dujella} are satisfied. Applying this 
	lemma to Inequality \eqref{Duj} with the following data
	
	\begin{align*}
		\tau:=\cfrac{\log  g }{\log \alpha_{\mathcal{N}}},\quad 
		\mu:=\cfrac{\log\biggl(\cfrac{ ( g -1)^3 a_{\mathcal{N}}}{d_1d_2d_3 
				\cdot  ( g^\ell-1   ) }\biggr)}{\log 
			\alpha_{\mathcal{N}}},\quad 
		A:=\cfrac{8}{ \log \alpha_{\mathcal{N}}},\quad  B:= g ,\quad 
		\text{and}\quad w:=m 
	\end{align*}
	with $\quad 1 \leq d_1  \leq d_2  \leq d_3 \leq  g -1$ and $1 \leq \ell 
	\leq 194$.
	\\
	As $k < 8 n \log g <1.99\times 10^{54}$, we can take $M:=1.99 \times 
	10^{54}$. With \textit{Mathematica} we get the following results :
	\begin{table}[ht]
		\caption{Upper bound on $m$}
		\label{ bound_m}
		\begin{center}
			\begin{tabular}{|c|c|c|c|c|c|c|c|c|c|}
				\hline \rule{0pt}{1.5 \normalbaselineskip}
				$ g $ &  {2} &   {3} &   {4} &   {5} &   {6} &   {7} &
				{8} &   {9} &   {10}
				\\ 
				
				\hline
				\hline \rule{0pt}{1.5 \normalbaselineskip}
				$q_t$  & $q_{118}$ & $q_{100}$ & $q_{110}$ & $q_{115}$ 
				& $q_{90}$ & $q_{106}$ & $q_{112}$ & $q_{102}$   &
				$q_{96}$
				\\ 
				\hline \rule{0pt}{1.5 \normalbaselineskip} 
				$\varepsilon  \ge$  &0.004 & 0.0007 & 0.0003 & 0.001
				& 0.0002 & 0.0005 & 0.0001 & 0.005 &0.001
				\\
				
				\hline \rule{0pt}{1.5 \normalbaselineskip}
				$m  \leq$  & 199  & 125    & 102    & 88       & 78  
				& 72     & 68      & 62   & 60    
				\\ 
				\hline
			\end{tabular}
		\end{center}
	\end{table}
	In all cases, we can conclude that 
	\[ 1 \leq m \leq \cfrac{\log((8/\log \alpha_{\mathcal{N}})\cdot q_{115}/0.0009)}
	{\log 2 }  \leq 200.
	\]
	\subsubsection*{Step~3.} Finally, to reduce the bound on $n$ we have to choose
	\[
	\Lambda_3:=\log (\Gamma_3+1)=(n)\log g -k \log \alpha_{\mathcal{N}}+\log 
	\biggl(\cfrac{ ( g -1)^3 a_{\mathcal{N}}}{d_1d_2d_3 \cdot   ( g^\ell-1   )
		(g^m-1   ) }\biggr).
	\]
	We have, 
	$$
	\biggl|  e^{-\Lambda_3}-1 \biggr|< \cfrac{2}{ g ^{n-1}}<\cfrac{1}{2},
	$$
	which is valid for $n \geq 4$ and $ g  \geq2$. It follows that 
	\begin{align}\label{Duj3}
		\Biggl| m   \cfrac{\log g}{\log \alpha_{\mathcal{N}}} - 
		k  + \cfrac{\log \biggl(\cfrac{ ( g -1)^3 a_{\mathcal{N}}}{d_1d_2d_3 
				\cdot  ( g^\ell-1   )   (g^m-1   ) }
			\biggr)}{\log \alpha_{\mathcal{N}}} 
		\Biggr|<\cfrac{4}{ g ^{n-1} \alpha_{\mathcal{N}} }. 
	\end{align}
	Now we have to apply Lemma~\ref{Dujella} to \eqref{Duj3} by taking the
	following
	parameters 
	\begin{align*}
		& \tau:=\cfrac{\log  g }{\log \alpha_{\mathcal{N}}},\quad 
		\mu:=\cfrac{\log \biggl(\cfrac{ ( g -1)^3 a_{\mathcal{N}}}{d_1d_2d_3 
				\cdot  (g^\ell-1   )   ( g^m-1   ) }
			\biggr)} {\log \alpha_{\mathcal{N}}}, 
		\quad A:=\cfrac{8}{ \log \alpha_{\mathcal{N}}},\quad 
		\\
		& B:= g ,\quad \text{and}\quad w:=n-1
	\end{align*}
	with $\quad 1 \leq d_1  \leq d_2  \leq d_3 \leq  g -1$, $ \quad 1 \leq 
	\ell  \leq 194$ and $1 \leq m \leq 183$. 
	
	Using the fact that $k < 8 n \log g <1.99\times 10^{54}$, we can take 
	$M:=1.99 \times 10^{54}$, and we get
	
	\begin{table}[ht] 
		\caption{Upper bound on $n$}
		\label{Bound_n}
		\begin{center}
			\begin{tabular}{|c|c|c|c|c|c|c|c|c|c|}
				
				\hline \rule{0pt}{1.5 \normalbaselineskip}
				{$ g $} &  {2} &   {3} &   {4} &   {5} &   {6} 
				&   {7} &   {8} &   {9} &   {10} 
				\\
				\hline
				\hline \rule{0pt}{1.5 \normalbaselineskip}
				{$q_t$}  & $q_{118}$ & $q_{99}$ & $q_{110}$ & 
				$q_{115}$ & $q_{90}$ & $q_{106}$ & $q_{112}$ & 
				$q_{102}$   & $q_{96}$
				\\
				\hline \rule{0pt}{1.5 \normalbaselineskip} 
				{$\varepsilon  \geq$} &0.00006 & 0.002 & 0.007 & 
				0.0002& 0.0008 & 0.002 & 0.0005 & 0.02 &0.009
				\\
				\hline \rule{0pt}{1.5 \normalbaselineskip}
				{$n  \leq$ } & 204  & 124   & 99   & 89      & 77
				& 71     & 67     & 61   & 59  
				\\  
				\hline
			\end{tabular}
		\end{center}
	\end{table}
	\vspace{3mm}
	It follows from the above table that
	$$
	1 \leq n \leq \cfrac{\log((4/\log \alpha_{\mathcal{N}})\cdot
		q_{118}/10^{-6})}{\log 2} \leq 205,
	$$ 
	which is valid for all $ g $ such that $2 \leq  g  \leq 10$. In light of the 
	above results, we need to check the equation \eqref{Meq1} in the cases $2
	\leq  g  \leq 10$ for $1  \leq d_1,d_2,d_3  \leq 9$, $1   \leq n   \leq 205$,
	$1   \leq m   \leq 200$, $1  \leq \ell  \leq 194$ and $1  \leq k  
	\leq 11500$. A quick inspection using \textit{Sagemath} reveals 
	that the Diophantine equation \eqref{Meq1} in the range $2  \leq g 
	\leq 10$ has only the solution listed in the statement of 
	Theorem~\ref{MT11}. This completes the proof of Theorem~\ref{MT11}.

	\section{Discussions}
	In addition to Baker's method and linear forms in logarithms, there are 
	other so-called ``classical'' methods and techniques for solving exponential
	Diophantine equations. These include the modular arithmetic method, 
	$p$-adic analysis, Fermat's method of infinite descent, the factorization 
	method, solving using inequalities, the mathematical induction method, 
	the parametric method, and so on. It would be interesting to treat the 
	same problems approached in this article with other methods than those 
	of the linear forms in logarithms. The modular arithmetic method could 
	be used to determine Narayana numbers, which are products of three 
	repdigits in base $g$ with $g\geq 2$ due to the interesting divisibility 
	properties possessed by the repdigits. 
	
	\section{ Acknowledgments}
	
	The authors express their gratitude to the anonymous reviewer for the
	instructive suggestions. The first author is partially supported by 
	Universit\'e de Kara (Togo) and the second author is supported by UCAD, 
	the Universit\'e Cheikh Anta Diop de Dakar. This project was initiated
	when the first author visited UCAD on a research stay. He thanks the 
	authorities for the warm hospitality and the working environment.

	\section*{Addresses}
	\noindent 
	$ ^{1} $ Facult\'e des Sciences et Techniques (FST), D\'epartement de Math\'ematiques, Universi\'e de Kara, Togo.
	\noindent
	$ ^{2} $ Laboratoire d’Algèbre, de Cryptologie, de G\'eom\'etrie Alg\'ebrique et Applications (LACGAA), Universit\'e Cheikh Anta Diop de Dakar (UCAD), S\'en\'egal.\\
	Email: \href{mailto:pagdame.tiebekabe@ucad.edu.sn}{pagdame.tiebekabe@ucad.edu.sn}
	\vspace{0.5cm}
	
	\noindent 
	$ ^{2} $  Laboratoire d’Algèbre, de Cryptologie, de G\'eom\'etrie Alg\'ebrique et Applications (LACGAA), Universit\'e Cheikh Anta Diop de Dakar (UCAD), S\'en\'egal
	
	\noindent
	Email: \href{mailto:kossirichmond.kakanou@ucad.edu.sn}{kossirichmond.kakanou@ucad.edu.sn}
	
	\vspace{0.5cm}
	
	\noindent $ ^{3} $ Faculty of Sciences Dhar El Mahraz, Sidi Mohamed Ben Abdellah University,
	Atlas-Fez, Morocco 
	
	\noindent	Email: \href{mailto:beyakouhamid@gmail.com}{beyakouhamid@gmail.com}
\end{document}